\documentclass[11pt]{amsproc}		
\usepackage{fancyhdr}
\usepackage{latexsym,amsfonts,amssymb,amsthm,amsmath}
\usepackage{stmaryrd} 
\usepackage{color}
\usepackage{pifont}
\usepackage[alphabetic]{amsrefs}

\usepackage[top=1.35in,bottom=1.35in,left=1.6in,right=1.6in]{geometry}
\usepackage[all]{xy}
\usepackage{todonotes}

\fancypagestyle{firststyle}{
	
	\fancyhf{}
	\fancyhead[C]{\footnotesize \textit{Enchiridion: Mathematical User's Guides} \hfill Volume 2 (September 2016)\\ \url{mathusersguides.com} \hfill \copyright 2016 CC BY-NC 4.0}
}

\newtheorem{lemma}{Lemma}[section]

\newtheorem{idea}[lemma]{Key Idea}  
\newtheorem{principle}[lemma]{Principle}
\newtheorem{ex}[lemma]{Example}

\newcommand{\F}{\mathbb{F}}
\newcommand{\Q}{\mathbb{Q}}
\newcommand{\R}{\mathbb{R}}
\newcommand{\Z}{\mathbb{Z}}

\newcommand{\Ainf}{\mathbb{A}_{\infty}}
\newcommand{\Einf}{\mathbb{E}_{\infty}}
\newcommand{\Hinf}{\mathbb{H}_{\infty}}

\newcommand{\de}{\delta}
\newcommand{\ep}{\epsilon}

\newcommand{\la}{\lambda}
\newcommand{\phy}{\varphi}
\newcommand{\Si}{\Sigma}
\newcommand{\te}{\theta}

\newcommand{\mm}{\mathfrak m} 

\newcommand{\llb}{\llbracket}
\newcommand{\op}{\oplus}
\newcommand{\ot}{\otimes}
\newcommand{\rrb}{\rrbracket}

\newcommand{\Alg}{\mathrm{Alg}}
\newcommand{\Alghat}{\widehat{\mathrm{Alg}}}
\newcommand{\Mod}{\mathrm{Mod}}
\newcommand{\Modhat}{\widehat{\mathrm{Mod}}}

\newcommand{\T}{\mathbb{T}}
\newcommand{\That}{\widehat{\mathbb{T}}}
\newcommand{\PP}{\mathbb{P}}
\newcommand{\Phat}{\widehat{\mathbb{P}}}

\newcommand{\ral}{\xrightarrow} 
\newcommand{\surj}{\twoheadrightarrow}

\newcommand{\id}{\mathrm{id}}

\newcommand{\abs}[1]{\lvert #1 \rvert}

\DeclareMathOperator{\coker}{coker}
\DeclareMathOperator{\colim}{colim}
\DeclareMathOperator{\Sym}{Sym}

\title[User's guide: Completed power operations]{A user's guide: Completed power operations for Morava $E$-theory}

\author{Martin Frankland}
\address{Institut f\"ur Mathematik\\
Universit\"at Osnabr\"uck\\
Osnabr\"uck, Germany}
\email{martin.frankland@uni-osnabrueck.de}

\begin{document}

\maketitle
\thispagestyle{firststyle}

\tableofcontents

\section{Key insights and central organizing principles}

This user's guide is about the paper \textit{Completed power operations for Morava $E$-theory}, written jointly with Tobias Barthel \cite{BarthelFrankland15}.

\subsection{Background} \label{sec:Background}

Given an $\Einf$ (or more generally $\Hinf$) ring spectrum $E$, \emph{power operations} for $E$ are the algebraic structure found in the homotopy of a commutative $E$-algebra. For $E = H\F_p$, the mod $p$ Eilenberg-MacLane spectrum, power operations generalize the Dyer--Lashof operations on the mod $p$ homology of infinite loop spaces. For $E = KU$, periodic complex $K$-theory, power operations on the homotopy of \emph{$p$-complete} $KU$-algebras have been described in work of McClure and Bousfield. They are given by $\te$-algebras over the $p$-adic integers $\Z_p$, which are commutative $\Z_p$-algebras with some additional structure.

From now on, fix a prime number $p$ and a height $h \geq 1$, and consider Morava $E$-theory $E = E_{h}$ at height $h$, and Morava $K$-theory $K(h)$. Power operations for $E$ on the homotopy of $K(h)$-local commutative $E$-algebras are understood thanks to work of Ando, Hopkins, Strickland, and Kashiwabara. Building on this, Rezk constructed a monad $\T$ on the category $\Mod_{E_*}$ of $E_*$-modules which encodes the power operations \cite{Rezk09}. For example, at height $h=1$, Morava $E$-theory is $p$-complete $K$-theory $E_1 = KU_p^{\wedge}$, and the monad $\T$ encodes $\te$-algebras.

Working $K(h)$-locally is crucial in those constructions. Now, an $E$-module spectrum $X$ is $K(h)$-local if and only if its homotopy $\pi_* X$ is an $E_*$-module which is $L$-complete, a property studied notably in work of Greenlees, May, Hovey, and Strickland\footnote{More precisely, $L$-completion is the best approximation of $\mm$-adic completion by a right exact functor.
The comparison map $L_0 M \to M^{\wedge}_{\mm}$ is an isomorphism for instance when the module $M$ is finitely generated or flat. In those cases, being $L$-complete is the same as being $\mm$-adically complete.}. This provides the insight motivating the paper:
\begin{quote}
\emph{If we only work with $K(h)$-local $E$-module spectra on the topological side, then we should only work with $L$-complete $E_*$-modules on the algebraic side.}
\end{quote}
If $A$ is a $K(h)$-local commutative $E$-algebra, then its homotopy $\pi_* A$ is an $L$-complete $E_*$-module which carries a $\T$-algebra structure. Since algebras over a monad have good categorical properties, it would be nice to encode the structure found in $\pi_* A$ via a monad. Namely, the monad $\T$ on $\Mod_{E_*}$ ought to be tightened to a monad $\That$ on the subcategory $\Modhat_{E_*}$ of $L$-complete $E_*$-modules, so that $\pi_* A$ is naturally a $\That$-algebra. The following diagram explains the problem schematically:
\[
\xymatrix @C=4.0pc @R=4.0pc{
\Alg_{E} \ar@<-0.7ex>[d]_{L_{K(h)}} & \Mod_{E_*} \ar@<-0.7ex>[d]_{L_0} \ar@(r,ur)_{\T} \\
**[l] \Alghat_{E} \ar@<-0.7ex>[u] \ar[ur]^{\pi_*} \ar[r]^-{\pi_*} & \Modhat_{E_*}. \ar@<-0.7ex>[u] \ar@(r,ur)_{\That} \\
}
\]

\begin{principle}
Monads are a convenient tool to encode additional structure, in topology as well as in algebra.
\end{principle}

In our situation, here is how monads encode the information on the topological side. Start with the free commutative $E$-algebra monad $\Sym \colon \Mod_{E} \to \Mod_{E}$ on the (model) category of $E$-modules, given by symmetric powers:
\[
\Sym(M) = \bigvee_{n \geq 0} \Sym^n M = \bigvee_{n \geq 0} \left( M^{\wedge_E n} \right)_{\Si_n}.
\]
Its derived functor $\PP \colon h\Mod_{E} \to h\Mod_{E}$ becomes a monad on the homotopy category of $E$-modules, namely the free $\Hinf$ $E$-algebra monad, given by \emph{extended powers}:
\[
\PP(M) = \bigvee_{n \geq 0} \PP_n M = \bigvee_{n \geq 0} \left( M^{\wedge_E n} \right)_{h \Si_n}.
\]
The functor $\T_n \colon \Mod_{E_*} \to \Mod_{E_*}$ is constructed so that there is an approximation map
\[
\T_n \left( \pi_* M \right) \to \pi_* L_{K(h)} \PP_n M
\]
for any $E$-module $M$, which is an isomorphism if $M$ is finitely generated and free, or \textbf{finite free} for short.

The following table summarizes the ingredients of the problem, and how the main theorem of \cite{BarthelFrankland15} improves the construction of the monad $\T$ to a monad $\That$ which takes $L$-completeness into account.

\bgroup
\def\arraystretch{1.6}
\begin{tabular}{|l||c|c|} \hline
& Topology & Algebra \\ \hline \hline
Ground ring & $E$ & $E_* := \pi_* E$ \\ \hline
Modules & $\Mod_{E}$ & $\Mod_{E_*}$ \\ \hline
Completion & $L_{K(h)} \colon \Mod_E \to \Modhat_E$ & $L_0 \colon \Mod_{E_*} \to \Modhat_{E_*}$ \\ \hline
Monad of interest & $\PP = \bigvee_{n \geq 0} \PP_n$ & $\T := \bigoplus_{n \geq 0} \T_n$ \\ \hline
Completed version & $\Phat := L_{K(h)} \PP$ & $\That = L_0 \T$ \\ \hline
Compatibility & $L_{K(h)} \PP \ral{\sim} L_{K(h)} \PP L_{K(h)}$ & $L_0 \T \ral{\text{iso?}} L_0 \T L_0$ \\
& & \textbf{Main Theorem: Yes!} \\ \hline
Consequence & $\Phat$ is a monad on $h\Modhat_E$. & $\That$ is a monad on $\Modhat_{E_*}$. \\
\hline
\end{tabular}
\egroup

There is a trade-off between algebraic invariants that encode more \emph{structure}, which are more powerful but hard to compute, and crude invariants, which are less powerful but easier to compute. However, there is no downside to encoding more \emph{properties}. As an analogy, if $X$ is a space, then it is sometimes useful to view $\pi_2 X$ as a $\pi_1 X$-module, sometimes more convenient to view it as an abelian group. However, there would be no benefit in forgetting that $\pi_2 X$ is abelian and viewing it merely as a group.

\begin{principle}
Algebraic invariants of topological objects should encode as much structure as possible (or as convenient), and as many properties as possible.  
\end{principle}

In the case of a $K(h)$-local commutative $E$-algebra $A$, one could look at all power operations on $\pi_* A$, encoded by a $\T$-algebra structure, or the underlying commutative $E_*$-algebra, or the underlying $E_*$-module, which only depends on the underlying $E$-module of $A$. In other words, there are forgetful functors:
\[
\xymatrix{
\Alg_{\T} \ar[r] & \Alg_{E_*} \ar[r] & \Mod_{E_*}, \\ 
}
\]
along with $L$-complete analogues, where we require the underlying $E_*$-module to be $L$-complete. The underlying $E_*$-module of $\pi_* A$ is $L$-complete, and we want to remember that \emph{property}.

\subsection{Algebraic ideas}

To prove the main theorem, our strategy is twofold: reduce the problem about $E_*$-modules to finite free modules, then prove the desired properties about $E$-modules, on the topological side. Here are the ingredients on the algebraic side.

\begin{idea}
Reduce statements about $L$-complete $E_*$-modules to statements modulo the maximal ideal $\mm \subset E_*$. 
\end{idea}

This idea is implemented by a Nakayama-type lemma \cite{BarthelFrankland15}*{Lemma~A.8}, which provides a useful tool to detect isomorphisms between $L$-complete modules. Using this, one then works with vector spaces over $E_* / \mm$, though the functors involved, such as $E_* / \mm \otimes \T_n$ are not additive.

\begin{idea}
Reduce the problem to simpler $E_*$-modules and simpler $E_*$-module maps.
\end{idea}

This idea is implemented in \cite{BarthelFrankland15}*{\S 4.1} using the following tricks.

\begin{itemize}
\item Any cokernel can be written (canonically) as a \emph{reflexive} coequalizer. In particular, any module is a reflexive coequalizer of free modules.
\item A free module is (canonically) a filtered colimit of finite free modules. More generally, a module is \emph{flat} if and only if it is a filtered colimit of finite free modules.
\item For a direct sum of modules $M \op N$, there is a ``binomial formula''
\[
\T_n (M \op N) = \bigoplus_{i+j=n} \T_i(M) \ot \T_j(N).
\]
This follows from the coproduct preservation formula $\T(M \op N) = \T(M) \ot \T(N)$.
\item The most relevant $E_*$-module maps are those given by multiplication by a scalar $\nu \in E_*$. In fact, it suffices to consider $\nu \in E_*$ taken from the (finite) set of generators of the maximal ideal $\mm \subset E_*$. 
\end{itemize}

\begin{idea}
Let $F$ be a finite free $E_*$-module. Then an $E_*$-module map $\phy \colon F \to F$ is nilpotent modulo the maximal ideal $\mm \subset E_*$ if and only if $\phy^{-1} F$ is trivial modulo $\mm$, i.e., the following equality holds:
\[
E_* / \mm \ot \colim \left( F \ral{\phy} F \ral{\phy} F \to \ldots \right) = 0. 
\]
\end{idea}

\subsection{Topological ideas}

The translation into topology relies on the fact that taking homotopy $\pi_* \colon \Mod_{E} \to \Mod_{E_*}$ induces an equivalence $h\Mod_{E}^{\text{ff}} \cong \Mod_{E_*}^{\text{ff}}$ between the homotopy category of finite free $E$-modules and the category of finite free $E_*$-modules.

\begin{idea}
Once the problem has been reduced to finite free modules, express algebraic problems about $E_*$-modules topologically in terms of $E$-modules.
\end{idea}

This strategy is carried out in \cite{BarthelFrankland15}*{\S 4.2}. For this, we need topological analogues of some of the algebraic facts mentioned above.

\begin{itemize}
\item An $E$-module is flat if and only if it is a filtered colimit of finite free $E$-modules.
\item On an $E_*$-module of the form $F = \pi_* M$, the $E_*$-module map $\nu \colon F \to F$ given by multiplication by a scalar $\nu \in E_*$ is the effect on homotopy of the $E$-module map $\nu \wedge_E \id \colon M \to M$.
\item If $f \colon M \to M$ is an $E$-module map, then the algebraic construction 
$(\pi_*f)^{-1} (\pi_*M)$ is the homotopy of the mapping telescope
\[
f^{-1} M = \colim \left( M \ral{f} M \ral{f} M \to \ldots \right)
\]
since $\pi_*$ preserves filtered colimits.
\end{itemize}

\subsection{Side issues}

In \cite{Rezk09}, treating the $\Z$-grading required significant work and care. In \cite{BarthelFrankland15}, the grading does not play an important role. Most of the argument can be made over the ring $E_0$. Minor adaptations then yield the result over the graded ring $E_*$, using the fact that $E_*$ is $2$-periodic.

\section{Metaphors and imagery}

\subsection{Rings and modules}

The paper relies on the fact that Morava $E$-theory is a highly structured ring spectrum, namely, an $\Einf$ ring spectrum. An $\mathbb{E}_1$ (or equivalently, $\Ainf$) ring structure guarantees that categories of module spectra behave well homotopically, while an $\Einf$ ring structure guarantees that the smash product $\wedge_E$ of $E$-modules behaves well. Since the advent of models for spectra with a strictly associative smash product, notably $S$-modules \cite{EKMM97} and symmetric spectra \cite{HoveySS00}, we can do stable homotopy theory in a way that mimics algebra. As far as I understand, this was part of the impetus behind ``brave new rings''---an expression coined by Waldhausen \cite{Greenlees07}---and higher algebra.

The metaphor is this: \emph{Ring spectra and module spectra are like ordinary rings and modules... but trickier}. At least the formal manipulations are similar. In our paper, we deal with notions for modules such as free, finitely generated, finitely presented (or perfect), and flat, as well as module maps given by multiplication by a scalar, a scalar acting invertibly on a module, etc. For these notions, the analogy between topology and algebra works remarkably well. In fact, our extension of Lazard's flatness criterion to module spectra \cite{BarthelFrankland15}*{\S 2.2} was motivated by this analogy, and in turn strengthens it.

As useful as it may be, the analogy has its limitations. We can think of extended powers $\PP_n M = (M^{\wedge n})_{h\Sigma_n}$ as the analogue of symmetric powers in algebra $\Sym^n M = (M^{\ot n})_{\Sigma_n}$. However, the homotopical construction is richer. For the monad $\T \colon \Mod_{E_*} \to \Mod_{E_*}$ described in Section~\ref{sec:Background}, every $\T$-algebra has an underlying commutative $E_*$-algebra, and the decomposition $\T = \oplus_{n \geq 0} \T_n$ corresponds to the decomposition $\Sym = \oplus_{n \geq 0} \Sym^n$. There is a natural comparison map $\Sym^n M \to \T_n M$ which is rarely an isomorphism, as $\T_n$ encodes additional, more exciting structure. For instance, at height $h=1$, the monad $\T$ is related to $\te$-rings; see \cite{BarthelFrankland15}*{Theorem~6.14} for a precise statement. A $\te$-ring is a commutative ring $R$ equipped with a non-linear operation $\te \colon R \to R$ satisfying certain equations. The operation $\te$ recovers the Adams operation $\psi \colon R \to R$ via the formula $\psi(x) = x^p + p \te(x)$.

Since the paper is about compatibility of the functors $\T_n$ with $L$-completion, let us describe metaphors related to $L$-completion, and then come back to the functors $\T_n$.

\subsection{The $p$-adic topology and power series}

As mentioned previously, $L$-completion is closely related to $\mm$-adic completion, and they often agree. At height $h=1$, the maximal ideal is $\mm = (p) \subset E_* = \Z_p [u^{\pm}]$, where $\Z_p$ denotes the $p$-adic integers. Here, let us ignore the $2$-periodicity of $E_*$ and focus on the degree zero part $E_0 = \Z_p$. Thus, $L$-completion is a variant of $p$-adic completion $M^{\wedge}_p = \lim_k M / p^k M$.

In the $p$-adic topology, $p$ is viewed as small, and higher powers $p^k$ are even smaller. Of course, this is unrelated to the ``size'' of $p^k$ as an integer number. As an analogy, think of the behavior of polynomials or power series in $x$ as $x$ tends to $0$. In that situation, $x$ is viewed as a small quantity, and $x^{100}$ is much smaller. The behavior of the power series is dictated by the term with lowest exponent, for example:
\[
f(x) = \overbrace{x^5}^{\text{dominant term}} + \quad \overbrace{x^{6} + 10 x^{23} + 3 x^{40} + \ldots}^{\text{negligible}}
\] 
Replacing $x$ by $p$ everywhere, we can say the same about the elements of $\Z_p$, and more generally within an $L$-complete $\Z_p$-module.

It turns out that this analogy can be made rigorous. In \cite{Rezk13ana}, Rezk shows that $L$-completion is given by \emph{analytic $p$-completion}, in the following sense:
\[
L_0 M \cong \coker \left( M \llb x \rrb \ral{(x-p)} M \llb x \rrb \right).
\]
One can then show that in an $L$-complete abelian group $M = L_0 M$, every element $f \in M$ admits a ``Taylor expansion around $p$'', which is a certain power series $\sum_n c_k (x-p)^k \in M \llb x \rrb$ representing $f$. An analogous description holds at higher height, where the maximal ideal $\mm$ has $h$ generators instead of one. $L$-completion is then analytic completion with respect to $h$ variables.

\subsection{Controlling error terms}

In a given $L$-complete module $M$, $p$ is ``small'' and $p^k$ is ``close'' to $0$ if $k$ is large. Given a map of modules $f \colon M \to N$, one can ask how continuous $f$ is with respect to the $p$-adic topology. Now we are interested in the functor $\T_n \colon \Mod_{E_*} \to \Mod_{E_*}$ and we push the analogy one level up. Think of the quotient map $q \colon M \surj M / p^k M$ as being ``close'' to an isomorphism when $k$ is large. In our paper, the main step towards the main theorem is to show that the induced map
\begin{equation}\label{eq:IsoAfterTensoring}
\xymatrix @C=4pc {
E_* / \mm \ot \T_n (E_*) \ar[r]^-{E_* / \mm \ot \T_n(q)} & E_* / \mm \ot \T_n (E_* / \mm^k) \\
}
\end{equation}
is an isomorphism for $k$ large enough \cite{BarthelFrankland15}*{Proposition~4.1}.

Here is an analogy with first-year calculus. Like many people, I think of the $\ep\text{--}\de$ definition of a limit $\lim_{x \to a} f(x) = L$ as a challenge and response. The $\ep$ is the challenge, while $\de = \de(\ep)$ is my response to that challenge. I don't expect $f(x) = L$ to hold anywhere, but as long as I can bound the error term $\abs{f(x)-L}$ by $\ep$, then I'm done.

In Equation~\eqref{eq:IsoAfterTensoring}, the fact that we are tensoring outside by the residue field $E_* / \mm$ is akin to a given challenge $\ep$. We could have tensored outside by $E_* / \mm^{100}$ for a more difficult challenge, akin to taking a smaller $\ep$, but in this situation, it turns out that tensoring outside by $E_* / \mm$ is enough, by a Nakayama-type argument. Our response, akin to finding a small enough $\de$, is to find an exponent $k$ large enough such that the ``small'' map $E_* \to E_* / \mm^k$ is sent to the map $\T_n(E_*) \to \T_n(E_* / \mm^k)$ which is not quite an isomorphism, but ``small enough'' that it becomes an isomorphism after tensoring with $E_* / \mm$.

Note that if $\T_n$ commuted with $E_* / \mm \ot -$, then the exponent $k=1$ would work and we'd be done. The functor $\Sym^n$ is an example of functor that commutes with $E_* / \mm \ot -$. One can show that $\T_n$ agrees with $\Sym^n$ (on finitely generated modules) for $n < p$. Hence, the first interesting example, where an exponent $k > 1$ might be required, is $\T_p$.

\begin{ex}\label{ex:Height1}
At height $h=1$, the quotient map $q \colon \Z_p \surj \Z/p^k$ induces the map
\[
\xymatrix{
\T_p(\Z_p) \ar[d]_{\simeq} \ar@{->>}[r]^-{\T_p(q)} & \T_p(\Z/p^k) \ar[d]^{\simeq} \\
\Z_p \op \Z_p \ar@{->>}[r] & \Z/p^{k+1} \op \Z/p^{k-1} \\
}
\]
where the bottom map is the direct sum of the two quotient maps. Therefore, the map $\Z/p \ot \T_p(q)$ is an isomorphism for all $k \geq 2$. The following table illustrates the thought process in a slightly cartoonish way.

\bgroup
\def\arraystretch{2.0}
\begin{tabular}{|l|c|c|} \hline
$k$ & $\Z/p \ot \T_p(q)$ & My thought \\ \hline \hline
$k=1$ & $\Z/p \op \Z/p \surj \Z/p$ & \textcolor{red}{Danger! Not yet injective.} \\ \hline
$k=2$ & $\Z/p \op \Z/p \ral{\simeq} \Z/p \op \Z/p$ & \textcolor{green}{\checkmark Safe.} The map is now injective. \\ \hline
$k=\text{bajillion}$ & $\Z/p \op \Z/p \ral{\simeq} \Z/p \op \Z/p$ & \textcolor{green}{\checkmark Even safer.} \\ \hline
\hline
\end{tabular}
\egroup
\end{ex}

\section{Story of the development}

While at Urbana--Champaign, I received an FQRNT Postdoctoral Fellowship starting in the fall 2011, to work under the supervision of Charles Rezk. In~\cite{Rezk09}*{\S 1.6}, Charles had raised a question for future research, which he suggested to me as a problem to work on, as it was a good fit for my interests. This problem became the topic of~\cite{BarthelFrankland15}. At first, I had to learn some background material on chromatic homotopy theory, Morava $E$-theory, power operations, $\te$-rings, $\la$-rings, $L$-complete modules, and so on. Charles' guidance was immensely helpful throughout the process.

At the beginning, I focused on the case of height $h=1$, in which case the monad $\T$ admits an explicit description, and thus the problem can be tackled using explicit formulas. Already at height $1$, the Nakayama-type reduction \cite{BarthelFrankland15}*{Lemma~A.8} was one of the first steps. Here's the reason: The goal was to show that the map
\[
\T_n(\eta) \colon \T_n M \to \T_n L_0 M
\]
induces an isomorphism upon applying $L$-completion $L_0$. The reduction says that under mild assumptions, it suffices to check this modulo the maximal ideal $(p) \subset \Z_p$. It was also clear that controlling the torsion in $\T_n(\Z/p^k)$ was a key issue, and that this module is computed by turning the cokernel diagram
\[
\xymatrix{
\Z_p \ar[r]^-{p^k} & \Z_p \ar@{->>}[r] & \Z/p^k \\
}
\]
into a reflexive coequalizer.

At the ``Strings and Automorphic Forms in Topology'' conference in Bochum in August 2012, I presented the result for $h=1$ in a contributed talk. My talk caught the attention of Tobias Barthel, then a graduate student at Harvard, who had been thinking about related topics. Thus, we started discussing on that occasion. At the Quillen Memorial Conference at MIT in October 2012, I gave a similar talk in a discussion session. There, I had a chance to discuss the project with Mark Hovey, Mark Behrens, and Haynes Miller, and talk some more to Tobias. Later in October, at the Midwest Topology Seminar at Michigan State University, Tobias and I discussed how my work at height $h=1$, along with the use of flat modules, could reduce the problem at arbitrary height $h$ to a key technical property of the functors $\T_n$; this reduction step later became \cite{BarthelFrankland15}*{\S 4.1}. That is when we ``officially'' joined forces to solve the general case $h \geq 1$.

In parallel, we had been wondering why Charles' proof that $\T$ is a monad does not also show formally that $\That$ is a monad. Looking at his proof more closely, we pinpointed where the argument breaks down for $\That$, because of peculiar features of the category $\Modhat_{E_*}$ of $L$-complete $E_*$-modules. Specifically, $E_*$ is small as an object of $\Mod_{E_*}$ but \emph{not} as an object of $\Modhat_{E_*}$; c.f. \cite{BarthelFrankland15}*{Remark~3.21}. This spurred us to study $L$-complete modules more thoroughly. Many results had already been worked out in~\cite{HoveyS99}*{Appendix~A}. Another important ingredient was found in unpublished notes of Hovey: $L$-completion preserves flatness \cite{BarthelFrankland15}*{Proposition~A.15}. We collected useful facts about $L$-complete modules into~\cite{BarthelFrankland15}*{Appendix~A}.

Proving the key property of $\T_n$ involved some reverse engineering. Looking at how $\T_n(E_*/\mm^k)$ is computed as a reflexive coequalizer, we saw that the important ingredient was nilpotency of the map $\T_n(\nu) \colon \T_n M \to \T_n M$ modulo the maximal ideal $\mm \subset E_*$, for a scalar $\nu \in E_*$. This in turn could be proved by going back to the topological side and using a mapping telescope \cite{BarthelFrankland15}*{Corollary~4.5}. Charles' input also helped us in that section, notably in the proof of \cite{BarthelFrankland15}*{Lemma~4.4}. For instance, he reminded us that a scalar $\nu \in E_*$ acts invertibly on an $E$-module $M$ if and only if the natural map 
\[
M = E \wedge_{E} M \to (\nu^{-1} E) \wedge_{E} M = \nu^{-1} M
\]
is an equivalence.

Working with flat $E_*$-modules made us think about flat $E$-module spectra. Helpful email discussions with Charles in August 2013 led to Lazard's criterion for flat module spectra \cite{BarthelFrankland15}*{\S 2.2}, which we use in the proof of~\cite{BarthelFrankland15}*{Corollary~4.5}. Said discussions were motivated by the following. We were comparing some work of Charles on flat modules over the (periodic) ring spectrum $E$ with work of Jacob Lurie on flat modules over \emph{connective} ring spectra. We wanted to clarify the relationship between the two approaches, and to find a common generalization; c.f. \cite{BarthelFrankland15}*{Remark~2.11}.

In \cite{BarthelFrankland15}*{\S 6}, we describe the monad $\T$ at height $h=1$. This result is known to experts and is due to Jim McClure \cite{BMMS86}*{\S IX}. However, it is not obvious how the calculations of McClure translate into the description of $\T$, and how this relates to work of Bousfield on $\te$-rings \cite{Bousfield96}, which is why we spelled out the details in our paper. Much of that section results from helpful email discussions with Jim McClure, spread out from October 2013 to February 2014. The first round of exchanges made it into the first arXiv version, posted in November 2013. Subsequent exchanges yielded the main changes in version 2, which was posted in February 2014. Notably, it was Jim McClure who suggested viewing his operation $Q$ as an operation on the homotopy of ($p$-complete) $\Hinf$ $KU$-algebras rather than on the (completed) $KU$-homology of $\Hinf$ ring spectra \cite{BarthelFrankland15}*{Proposition~6.8}.
 
The paper was submitted to Algebraic \& Geometric Topology in March 2014. We received the referee report in September 2014, submitted the revisions in October, and received the final acceptance in November.

\section{Colloquial summary}

I will start with a gentle introduction to topology.  Then I will leave a black box around the specifics of our paper \cite{BarthelFrankland15} and focus only on one aspect: the word ``completed'' in the title.

\subsection{What is topology?}

\begin{figure}[h]
\centering
\includegraphics{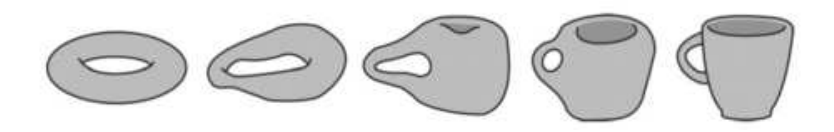}
\caption{A donut and a coffee mug are equivalent spaces.\newline
Image credit: Wikimedia Commons, via www.functionspace.com.}
\label{fig:DonutMug}
\end{figure}

Topology is the study of spaces, such as curves, surfaces, and higher-dimensional analogues. Unlike geometry, which cares about angles, lengths, and volume, topology looks at the qualitative features of a space, its general shape: how many pieces there are, whether there are holes, the number of holes, and so on. 
That's why topology is sometimes called rubber-sheet geometry: spaces can be stretched, compressed, bent, and it's all the same to us. According to a classic joke, a topologist cannot distinguish a donut from a coffee mug (Figure~\ref{fig:DonutMug}).

\begin{figure}[h]
\centering
\includegraphics[width=0.9\linewidth]{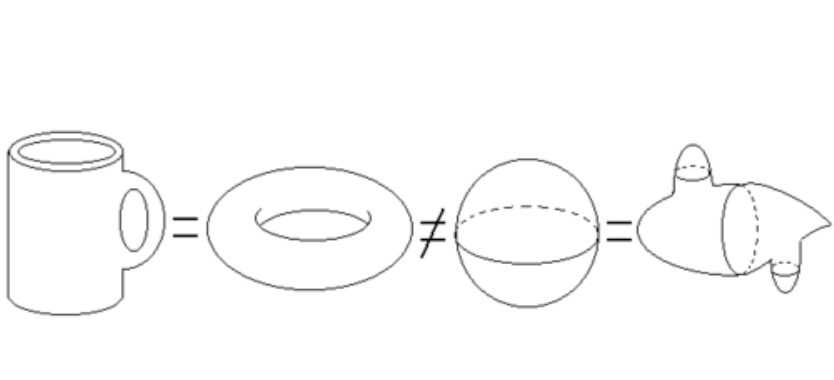}
\caption{A donut and a sphere are \emph{not} equivalent spaces.\newline
Image credit: www.functionspace.com.}
\label{fig:TorusSphere}
\end{figure}

Likewise, a balloon, the surface of a football, or the surface of the Earth are all the same space, namely a $2$-dimensional sphere.
However, the donut and the sphere are different spaces (Figure~\ref{fig:TorusSphere}).

\begin{figure}[h] 
\centering
\includegraphics{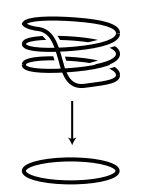}
\caption{A loop winding $3$ times around the circle.\newline
Image credit: Allen Hatcher, \textit{Algebraic Topology}, Cambridge University Press.}
\label{fig:CoveringCircle}
\end{figure}

Algebraic topology is the branch of topology that describes the shape of spaces using algebraic 
invariants, quantities that we can compute. For example, using algebra, we can make precise the idea that the donut has a hole in it that the sphere doesn't have. 
For another example, a loop in a circle 
is described by the number of times it's winding around the circle (Figure~\ref{fig:CoveringCircle}). 
By counting the winding number, we describe topological information (the loops in a circle) via an algebraic structure (the integer numbers $\Z$). 

As exciting as spaces are, you may wonder what they're good for.  Topology has become an important branch of mathematics, with connections to geometry, analysis, algebra, and mathematical physics.  Since the 1990s, algebraic topology has also been applied to problems such as data analysis, sensor networks, and robot motion planning.

\subsection{What is completeness?}

The paper \cite{BarthelFrankland15} takes place in a branch of algebraic topology called \emph{chromatic homotopy theory}. We study certain kinds of ``spaces'' and their algebraic invariants, which have the interesting feature of being \emph{complete}. Our main theorem says that certain computations can be made using only algebraic structures that are complete. 

In this section, I will sketch what ``complete'' means, and why it might be an interesting feature.

\begin{figure}[h]
\centering
\includegraphics[width=0.4\linewidth]{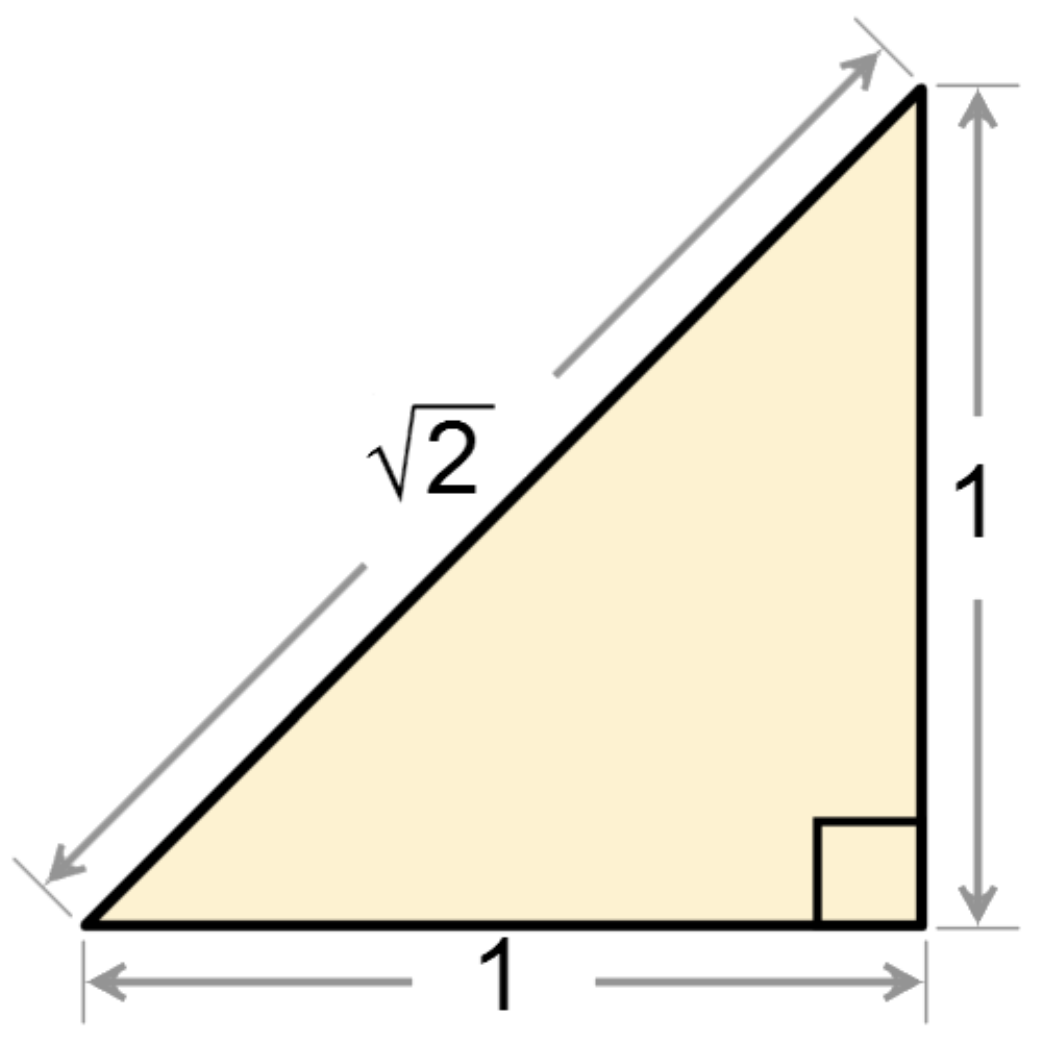}
\caption{A right triangle.\newline
Image credit: Wikimedia Commons.}
\label{fig:Triangle}
\end{figure}

Consider the number $\sqrt{2}$, which you may have seen in a geometry class (Figure~\ref{fig:Triangle}). 

Its decimal expansion is
\begin{align*}
\sqrt{2} &= 1.41421\ldots \\
&= 1 + 0.4 + 0.01 + 0.004 + 0.0002 + 0.00001 + \ldots \\
&= 1 \cdot 10^0 + 4 \cdot 10^{-1} + 1 \cdot 10^{-2} + 4 \cdot 10^{-3} + 2 \cdot 10^{-4} + 1 \cdot 10^{-5} + \ldots
\end{align*}
Recall that $\sqrt{2}$ is irrational, i.e., cannot be expressed as a fraction\footnote{This fact was discovered in Ancien Greece in the 5th century BC, sometimes attributed to Hippasus of Metapontum. According to legend, he was murdered for his discovery of irrational numbers, deemed shocking at the time by his fellow Pythagoreans.}. 
Another way to say this is that the decimal expansion of $\sqrt{2}$ does not start repeating after a while, which is what happens with a fraction, e.g.:
\[
\frac{5}{6} = 0.8\overline{3} = 0.833333333\ldots
\]
The rational numbers $\Q$ are \emph{incomplete}: there is a gap where $\sqrt{2}$ would be expected. 
In contrast, the real numbers $\R$ are \emph{complete}: they contain all numbers obtained from decimal expansions.

Another example of infinite sum appears in calculus, when expressing functions as their \emph{Taylor series}\footnote{Named after the English mathematician Brook Taylor, who introduced them formally in 1715. Some cases had already been used by Isaac Newton and James Gregory in the 1660s, and by Indian mathematician Madhava in the 14th century.}. 
For instance, the exponential function is
\[
e^x = 1 + x + \frac{x^2}{2!} + \frac{x^3}{3!} + \frac{x^4}{4!} + \frac{x^5}{5!} + \cdots
\]
This infinite sum is useful both in theory and in practice. To this day, some calculators use this formula to compute the exponential function, by adding the first few terms, depending how much precision we need. 
Take for instance $x=0.1$, and compare these values:
\begin{align*}
&e^{0.1} \approx 1.105171 \\
&1 + 0.1 + \frac{(0.1)^2}{2!} + \frac{(0.1)^3}{3!} \approx 1.105167.
\end{align*}
The first few terms already yield a good approximation.

In our paper, the word ``complete'' means that certain infinite sums are available, which is convenient. 
Let's say we compute something involving 
the algebraic invariants of our spaces. 
If we insist on a certain precision in the answer, we know we can achieve it by taking enough terms of the infinite sum.

Here's a non-mathematical analogy as to why completeness is useful. We topologists study spaces, as an ornithologist studies birds, or a botanist studies plants. The algebraic invariants we associate to a space provide some information about the space, as pictures of birds or plants provide some information about them. Having algebraic invariants that are complete would be like having high-resolution pictures, where one can zoom in and still see a clear picture.

\end{document}